\input amstex
  \documentstyle{amsppt}
   \magnification=\magstephalf
  \def \m{{\text {\bf m}}}
   
   \def \q{{\text {\bf q}}}
   
   \def \q{{ \text{\bf q}}}
   \def \p{{ \text{\bf p}}}
  
   \def \m{{ \text{\bf m}}}

   \def \N{{\Bbb N}}

\def \sC{{\Cal C}}

  \long\def\alert#1{\smallskip\line{\hskip\parindent\vrule
\vbox{\advance\hsize-2\parindent\hrule\smallskip\parindent.4\parindent
  \narrower\noindent#1\smallskip\hrule}\vrule\hfill}\smallskip}
   \def \hgt{\operatorname{ ht}}
   \def \dim{\operatorname{ dim}}
   \def \ker{\operatorname{ ker}}

   \def \Spec{\operatorname{ Spec}}

  \NoBlackBoxes
  \rightheadtext{ Generic fiber rings }\topmatter
\leftheadtext{Heinzer, Rotthaus and Wiegand } \topmatter

   \title
Generic fiber rings  of mixed \\
power series/polynomial rings
  \endtitle

   \author
   William Heinzer, Christel Rotthaus and Sylvia
   Wiegand
  \endauthor


   \thanks{The authors are grateful for
  the hospitality and cooperation of Michigan State, Nebraska
   and Purdue, where several work sessions on this research were
   conducted. Wiegand thanks the National Security Agency for its
support.}
   \endthanks

   \address{Department of Mathematics, Purdue
   University,
   West Lafayette, IN 47907-1395}
   \endaddress

   \address{Department of Mathematics, Michigan State University, East
  Lansing,
   MI 48824-1027}
   \endaddress

   \address{Department of Mathematics and Statistics,
   University of Nebraska,
   Lincoln, NE 68588-0323}
   \endaddress

\abstract{Let $K$ be a field, $m$ and $n$ positive integers, and
$X = \{x_1, \ldots, x_n\}$, and  $Y = \{y_1, \ldots, y_m\}$ sets of
independent variables over $K$. Let $A$ be the
localized polynomial ring $K[X]_{(X)}$. We
prove that every prime ideal $P$ in
$\widehat A = K[[X]]$
that is maximal with respect to $P \cap A = (0)$
has height $n-1$. We consider the mixed
power series/polynomial rings
$B := K[[X]]\,[Y]_{(X, Y)}$ and
$C := K[Y]_{(Y)}[[X]]$. For each prime
ideal $P$ of $\widehat B = \widehat C$
that is maximal with respect to either
$P \cap B = (0)$ or $P \cap C = (0)$, we
prove that $P$  has
height $n+m-2$. We also prove
each prime ideal   $P$ of $K[[X,Y]]$ that is maximal
with respect to $P \cap K[[X]] = (0)$ is of
height either $m$ or $n+m-2$. }
\endabstract

   \endtopmatter

   \document

\baselineskip 17pt

\subheading{1. Introduction and Background}

Let $(R,\m)$ be a Noetherian local integral domain and let
$\widehat R$ denote the $\m$-adic completion of $R$. The {\it
generic formal fiber ring} of  $R$ is the localization $(R
\setminus (0))^{-1}\widehat R$ of $\widehat R$. The {\it formal
fibers} of $R$ are the fibers of the morphism $\Spec \widehat R
\to \Spec R$; for a prime ideal $P$ of $R$, the formal fiber over $P$ is
$\Spec(\,(R_P/PR_P)\otimes_R\widehat R$\,). The formal fibers
 encode important information about the structure of $R$.
For example, the local ring $R$ is excellent provided it is
universally catenary and has geometrically regular
formal fibers \cite{2, (7.8.3), page 214}.

Let $R \hookrightarrow S$ be an injective homomorphism of
commutative rings. If $R$ is an integral domain, the {\it generic
fiber ring} of the map $R \hookrightarrow S$ is the localization
$(R \setminus (0))^{-1}S$ of $S$.
 In this article we study generic
fiber rings for  ``mixed" polynomial and power series rings
over a field. More precisely, for $K$ a field,  $m$ and $n$
positive integers, and $X=\{x_1, \ldots, x_n\}$ and $Y=\{y_1,
\ldots, y_m\}$ sets of variables over $K$, we consider the local
rings $A := K[X]_{(X)}$, $B := K[[X]]\,[Y]_{(X,Y)}$ and $C :=
K[Y]_{(Y)}[[X]]$, as well as their completions $\widehat A = K[[X]]$
and $\widehat B = \widehat C = K[[X, Y]]$.
Notice that there is a canonical inclusion map
$B \hookrightarrow C$.

\pagebreak

 We have the following local embeddings.

 $$\align
A &:= K[X]_{(X)}\hookrightarrow \widehat A:=K[[X]],\quad
\widehat A \hookrightarrow \widehat B = \widehat C = K[[X, Y]]
\quad \text{ and } \\
B &:= K[[X]]\,[Y]_{(X,Y)}\hookrightarrow C :=
K[Y]_{(Y)}[[X]]\hookrightarrow\widehat B=\widehat C=
K[[X]]\,[[Y]].
\endalign
$$

Matsumura proves in \cite{7} that the generic formal fiber ring of
$A$ has dimension $n-1 = \dim A - 1$,  and the generic formal
fiber rings of $B$ and $C$ have dimension  $n+m-2 = \dim B - 2 =
\dim C - 2$. However he does not address the question of whether
all maximal ideals of the generic formal fiber rings for $A$, $B$
and $C$ have the same height. If the field $K$ is countable, it
follows from \cite{3, Prop. 4.10, page 36} that all maximal ideals
of the generic formal fiber ring of $A$ have the same height.

In answer to a question raised by Matsumura in \cite{7},  Rotthaus
in \cite{10} establishes the following result. Let $n$ be a
positive integer. Then there exist excellent regular local rings
$R$ such that $\dim R = n$ and such that the generic formal fiber
ring of $R$ has dimension $t$, where the value of $t$ may be taken
to be any integer between $0$ and $\dim R - 1$. It is also shown
in \cite{10, Corollary 3.2} that there exists an excellent regular
local domain having the property that its generic formal fiber
ring contains maximal ideals of different heights.

Let $\widehat T$ be a complete Noetherian local ring and let $\sC$
be a finite set of incomparable prime ideals of $\widehat T$.
Charters and Loepp in \cite{1} (see also \cite{6, Theorem 17})
determine necessary and sufficient conditions for $\widehat T$ to
be the completion of a Noetherian local domain $T$ such that the
generic formal fiber of $T$ has as maximal elements precisely the
prime ideals in $\sC$. If $\widehat T$ is of characteristic zero,
Charters and Loepp give necessary and sufficient conditions to
obtain such a domain $T$ that is excellent. The finite set $\sC$
may be chosen to contain prime ideals of different heights. This
provides many examples where the generic formal fiber ring
contains maximal ideals of different heights.

Our main results may be summarized as follows.

\proclaim{1.1 Theorem} With the above notation, we prove that all
maximal ideals  of the generic formal fiber rings of $A$, $B$ and
$C$ have the same height. In particular, we prove: \roster \item
If $P$ is a  prime ideal of $\widehat A$ maximal with respect to
$P\cap A=(0),$ then $\hgt(P)=n-1$. \item If $P$ is a  prime ideal
of $\widehat B$ maximal with respect to $P\cap B=(0),$ then
$\hgt(P)=n+m-2$. \item If $P$ is a  prime ideal of $\widehat C$
maximal with respect to $P\cap C=(0),$ then $\hgt(P)=n+m-2$.
\item In addition, there are at most two possible values for the
height of a maximal 
ideal of  the generic fiber ring $(\widehat A
\setminus (0))^{-1}\widehat C$ of the inclusion map $\widehat A
\hookrightarrow \widehat C$. 

$\phantom {x}${\rm (a)} If $n\ge 2$ and
$P$ is a prime ideal of $\widehat C$ maximal with respect to \newline 
$\phantom {xxxxxx} P
\cap \widehat A = (0)$, then either $\hgt P = n+m-2$ or $\hgt P =
m$. 

$\phantom {x}${\rm (b)}  If $n = 1$, then all maximal ideals of the generic
fiber ring \newline
 $\phantom {xxxxxx}(\widehat A \setminus (0))^{-1}\widehat C$ have
height~$m$.
\endroster\endproclaim

We were motivated  to consider generic  fiber rings for the
embeddings displayed above  because of questions related to \cite
{4} and \cite{5} and ultimately because of the following  question
posed by Melvin Hochster.

\subheading{1.2 Question} Let $R$ be a complete local domain.
Can one describe or somehow classify the local maps of $R$ to
a complete local domain $S$ such that $U^{-1}S$  is a field,
where $U = R \setminus (0)$, i.e., such that the generic fiber
of $R \hookrightarrow S$ is trivial?

Hochster remarks  that if,  for example,
$R$  is equal characteristic  zero, one obtains such
extensions by starting with
$$
R = K[[x_1, ...,  x_n]] \hookrightarrow T =
L[[x_1, ..., x_n, y_1, ..., y_m]] \to T/P = S,
\tag{1.2.1}
$$
where   $K$  is a subfield of
$L$,  the  $x_i, y_j$ are formal indeterminates,  and  $P$ is a
prime ideal of  $T$  maximal with respect to being
disjoint from the image of   $R\setminus \{0\}$.
Of course, such prime ideals $P$
correspond to the maximal ideals of the generic fiber
$(R \setminus (0))^{-1} T$.

In Theorem 7.2, we answer  Question 1.2 in the special case where the extension
arises from the embedding in (1.2.1) with the field $L=K.$
We prove in this case that the dimension of the extension ring
$S$ must be either $2$ or $n$.

In \cite{5} we study extensions of integral domains
$R\overset{\varphi}\to\hookrightarrow S$ such that, for every
nonzero $Q  \in \Spec S$, we have $Q\cap R\not=(0).$ Such
extensions are called {\it trivial generic fiber extensions} or
{\it TGF extensions} in \cite{5}. One obtains such an extension by
considering a composition $R\hookrightarrow T \to T/P=S,$ where
$T$ is an extension ring of $R$  and $P \in \Spec T$ is maximal
with respect to $P\cap R=(0)$. Thus the generic fiber ring and so
also Theorem 1.1 give information regarding TGF extensions in the
case where the smaller ring is a mixed polynomial/power series
ring.

In addition,  Theorem 1.1  is  useful in the study of (1.2.1),
because the map in (1.2.1) factors through:
$$
R= K[[x_1,\hdots,x_n]] \hookrightarrow
K[[x_1,\hdots,x_n]]\,[y_1,\hdots,y_m]\hookrightarrow T =
L[[x_1,\hdots,x_n,y_1,\hdots,y_n]].
$$

Section 2 contains implications of Weierstrass' Preparation
Theorem to the prime ideals of power series rings. We first prove
a   technical proposition regarding a change of variables that
provides a ``nice" generating set for a given prime ideal $P$ of a
power series ring; then in Theorem 2.3 we prove that, in certain
circumstances, a larger prime ideal can be found with the same
contraction as $P$ to a certain subring. In Sections 3 and 4, we
prove parts 2 and 3 of Theorem 1.1 stated above. In Section 5 we
use a result of Valabrega for the two-dimensional case. We then
apply this result  in Section 6 to prove part 1 of Theorem 1.1,
and in Section 7  we prove part 4.

\subheading{2. Variations on a theme of Weierstrass}

In this section, we apply the Weierstrass Preparation Theorem
\cite{12, Theorem 5, page 139, and Corollary 1, page 145} to
examine the structure of a given prime ideal $P$ in the power
series ring $\widehat A = K[[X]]$, where $X = \{x_1, \ldots ,
x_n\}$ is a set of $n$ variables over the field $K$.  Here $A =
K[X]_{(X)}$ is the localized polynomial ring in these variables.
Our procedure is to make a change of variables that yields a
regular sequence in $P$ of a nice form.

\subheading{2.1  Notation} By a {\it change of variables}, we  mean
a finite sequence of  `polynomial' change of variables of the type
described below, where $X=\{x_1,\hdots,x_n\}$ is a set of $n$
variables over the field $K$. For example, with $e_i, f_i \in \N$,
consider
$$\multline
\phantom{xxxx}x_1 \mapsto x_1 +x_n^{e_1}=z_1, \qquad x_2 \mapsto
x_2 + x_n^{e_2}=z_2,\qquad  \ldots ,\\ x_{n-1} \mapsto x_{n-1} +
x_n^{e_{n-1}}=z_{n-1}, \qquad  x_n \mapsto x_n =z_n,\phantom{xxxx}
\endmultline
$$ followed by:
$$
\multline \phantom{xxxxx}z_1 \mapsto z_1=t_1,\qquad z_2 \mapsto
z_2 + z_1^{f_2}=t_2,\qquad \ldots ,\\ z_{n-1} \mapsto z_{n-1} +
z_1^{f_{n-1}}=t_{n-1}, \qquad  x_n \mapsto z_n+ z_1^{f_n}
=t_n.\phantom{xxxxx}
\endmultline
 $$
Thus a change of variables defines an automorphism of $\widehat A$
that restricts to an automorphism of $A$.

We also consider a change of variables for subrings of $A$ and
$\widehat A$. For example, if $A_1= K[x_2,\hdots,x_n]\subseteq A$
and  $S =K[[x_2,\hdots,x_n]]\subseteq \widehat A$, then by a {\it
change of variables inside}  $A_1$ and $S$, we mean a finite
sequence of automorphisms of  $A$ and $\widehat A$ of the type
described above on $x_2, \ldots, x_n$ that leave the variable $x_1$
fixed. In this case we obtain an automorphism  of $\widehat A$ that
restricts to an automorphism on each of $S$, $A$ and $A_1$.

\proclaim{2.2 Proposition } Let $\widehat
A:=K[[X]]=K[[x_1,\hdots,x_n]]$ and let $P \in \Spec \widehat A$
with $x_1 \not\in P$ and $\hgt P = r$, where $1 \le  r \le n-1$.
There exists  a change of variables $x_1 \mapsto
z_1:=x_1$ ($x_1$ is fixed), $x_2 \mapsto z_2,\,\hdots, \,x_n
\mapsto z_n$ and a regular sequence  $f_1,\hdots,f_r\in P$ so that, upon
setting $Z_1 = \{z_1, \ldots, z_{n-r}\}$, $Z_2 = \{z_{n-r+1},
\ldots, z_n\}$ and $Z = Z_1 \cup Z_2$, we have
$$\alignat2
f_1&\in K[[Z_1]]\,[z_{n-r+1},\hdots,z_{n-1}]\,[z_n]\quad&&\text{is
monic as a polynomial in
$z_n$}\\
f_2&\in K[[Z_1]]\,[z_{n-r+1},\hdots,z_{n-2}]\,[z_{n-1}]
\quad&&\text{is monic as a polynomial
in $z_{n-1}$, etc}\\
 &\vdots && \\
f_r&\in K[[Z_1]]\,[z_{n-r+1}] \quad&&\text{is monic as a
polynomial in $z_{n-r+1}$.}
\endalignat
$$
In addition:
 \roster
\item $P$ is a minimal prime of the ideal $(f_1,\hdots,f_r)\widehat A$.
\item  The $(Z_2)$-adic completion  of $K[[Z_1]]\,[Z_2]_{(Z)}$
is identical to the $(f_1,\hdots,f_r)$-adic completion
 and both equal  $\widehat A=K[[X]] = K[[Z]]$.
\item  If
$P_1 := P \cap K[[Z_1]]\,[Z_2]_{(Z)}$, then $P_1\widehat A = P$,
that is, $P$ is extended from $K[[Z_1]]\,[Z_2]_{(Z)}$.
\item The ring extension:
$$ K[[Z_1]]\hookrightarrow K[[Z_1]]\,[Z_2]_{(Z)}/P_1\cong K[[Z]]/P
$$
is finite (and integral).
\endroster
\endproclaim

\demo{Proof}  Since $\widehat A$ is a unique factorization domain,
there exists  a nonzero prime element $f$ in $P$. The power series
$f$ is therefore not a multiple of $x_1$, and so $f$ must contain
a monomial term $x_2^{i_2}\hdots x_n^{i_n}$ with a nonzero
coefficient in $K$. This nonzero coefficient in $K$ may be assumed
to be 1. There exists an automorphism $\sigma : \widehat A \to
\widehat A$ defined by the change of variables:
$$ x_1 \mapsto x_1 \qquad x_2 \mapsto t_2:=x_2 +
x_n^{e_2}\quad \hdots \quad x_{n-1} \mapsto t_{n-1}:=x_{n-1} +
x_n^{e_{n-1}}\qquad x_n \mapsto x_n,
$$
with $e_2,\hdots,e_{n-1}\in \Bbb N$ chosen suitably so that $f$
written as a power series in the variables
$x_1,t_2,\hdots,t_{n-1},x_n$ contains a term $a_nx_n^{s_n}$, where
 $s_n$ is a positive integer, and $a_n\in K$ is nonzero.
We assume that the integer $s_n$ is minimal among all integers $i$
such that a term $ax_n^i$ occurs in $f$ with a nonzero coefficient $a\in K$;  we
further assume that the coefficient $a_n=1$. By Weierstrass we
have that:
$$ f = m \epsilon,$$
where $m\in K[[x_1,t_2,\hdots,t_{n-1}]]\,[x_n]$ is a monic
polynomial in $x_n$ of degree $s_n$ and $\epsilon$ is a unit in
$\widehat A$. Since $f \in P$ is a prime element,  $m \in P$ is also a prime
element.  Using
Weierstrass again, every element $g\in P$ can be written as:
$$ g = mh + q,$$
where $h\in K[[x_1,t_2,\hdots,t_{n-1},x_n]] = \widehat A$ and
$q\in K[[x_1,t_2,\hdots,t_{n-1}]]\,[x_n]$ is a polynomial in $x_n$
of degree less than $s_n$. Note that
$$K[[x_1,t_2,\hdots,t_{n-1}]]\hookrightarrow
K[[x_1,t_2,\hdots,t_{n-1}]]\,[x_n]/(m)$$ is an integral (finite)
extension. Thus the ring $K[[x_1,t_2,\hdots,t_{n-1}]]\,[x_n]/(m)$
is complete. Moreover, the two ideals $(x_1,t_2,\hdots,t_{n-1},m)$
= $(x_1,t_2,\hdots,t_{n-1},x_n^{s_n})$ and
$(x_1,t_2,\hdots,t_{n-1},x_n)$ of $B_0 :=
K[[x_1,t_2,\hdots,t_{n-1}]]\,[x_n]$ have the same radical. Therefore
$\widehat A$ is the $(m)$-adic and the $(x_n)$-adic completion of
$B_0$ and  $P$ is extended from $B_0$.

This implies the statement for $r=1$, with $f_1=m, \,z_n=x_n,\,
z_1=x_1,\, z_2=t_2,\,\hdots,\,z_{n-1}=t_{n-1}$, \,$Z_1=\{
x_1,t_2,\hdots,t_{n-1}\}$ and $Z_2=\{z_n\}=\{x_n\}$. In
particular, when $r=1$, $P$ is minimal over $m\widehat A$,
so $P = m\widehat A$.

For $r>1$ we continue by induction on $r$. Let $P_0:=P\cap
K[[x_1,t_2,\hdots, t_{n-1}]]$. Since $m\notin K[[x_1,t_2,\hdots,
t_{n-1}]]$ and $P$ is extended from $B_0:=K[[x_1,t_2,\hdots,
t_{n-1}]]\,[x_n]$, then $P\cap B_0$ has height $r$ and $\hgt
P_0=r-1$. Since $x_1 \notin P$, we have $x_1\notin P_0$, and by
the induction hypothesis there is a change of variables $t_2
\mapsto z_2, \hdots , t_{n-1}\mapsto z_{n-1}$ of
$K[[x_1,t_2,\hdots,t_{n-1}]]$ and elements $f_2,\hdots,f_r\in P_0$
so that:

$$\alignat2
f_2&\in
K[[x_1,z_2\hdots,z_{n-r}]]\,[z_{n-r+1},\hdots,z_{n-2}]\,[z_{n-1}]\quad&&\text{is
monic in
$z_{n-1}$}\\
f_3&\in K[[x_1,z_2\hdots,z_{n-r}]]\,[z_{n-r+1},\hdots,z_{n-3}]\,[z_{n-2}]
\quad&&\text{is monic in $z_{n-2}$, etc}\\
 &\vdots && \\
f_r&\in K[[x_1,z_2,\hdots,z_{n-r}]]\,[z_{n-r+1}] \quad&&\text{is
monic in $z_{n-r+1}$,}
\endalignat
$$
and $f_2,\ldots,f_r$ satisfy the assertions of Proposition 2.2 for
$P_0$.

It follows  that $m,f_2,\hdots,f_r$ is a regular sequence of
length $r$ and that $P$ is a minimal prime of the ideal
$(m,f_2,\hdots,f_r)\widehat A$. Set $z_n= x_n$. We now prove  that
$m$ may be replaced by a polynomial $f_1\in
K[[x_1,z_2,\hdots,z_{n-r}]]\,[z_{n-r+1},\hdots,z_n]$. Write
$$ m = \sum_{i=0}^{s_{n}} a_iz_{n},$$
where the  $a_i\in K[[x_1,z_2,\hdots,z_{n-1}]]$. For each $i<s_n$,
apply  Weierstrass to $a_i$ and $f_2$ in order to obtain:
$$ a_i = f_2 h_i + q_i,$$
where $h_i$ is a power series in $K[[x_1,z_2,\hdots,z_{n-1}]]$ and
$q_i\in K[[x_1,z_2,\hdots,z_{n-2}]]\,[z_{n-1}]$ is a polynomial in
$z_{n-1}$. With $q_{s_n} = 1 =a_{s_n}$,  we define
$$ m_1 = \sum_{i=0}^{s_{n}} q_iz_{n}^i.$$
Now $(m_1,f_2,\hdots,f_r)\widehat A =(m,f_2,\hdots,f_r)\widehat A$
and we may replace $m$ by $m_1$ which is a polynomial in $z_{n-1}$
and $z_n$. To continue, for each $i<s_n$, write:

$$ q_i = \sum_{j,k} b_{ij} z_{n-1}^j \quad
\text{with}\quad b_{ij}\in K[[x_1,z_2,\hdots,z_{n-2}]].$$ For each
$b_{ij}$, we apply Weierstrass to  $b_{ij}$ and $f_3$ to obtain:
$$ b_{ij} = f_3h_{ij} + q_{ij},$$
where $q_{ij}\in K[[x_1,z_2,\hdots,z_{n-3}]]\,[z_{n-2}]$. Set
$$m_2 = \sum_{i,j} q_{ij} z_{n-1}^jz_n^i \in
K[[x_1,z_2,\hdots,z_{n-3}]]\,[z_{n-2},z_{n-1},z_n]$$ with $q_{s_n0}
=1$. It follows that $(m_2,f_2,\hdots,f_r)\widehat A
=(m,f_2,\hdots,f_r)\widehat A$. Continuing this process by
applying Weierstrass to the coefficients of
$z_{n-2}^kz_{n-1}^jz_n^i$ and $f_4$, we establish the existence of
a polynomial  $f_1\in K[[Z_1]]\,[z_{n-r+1},\hdots,z_n]$ that is
monic in $z_n$ so that $(f_1,f_2,\hdots,f_r)\widehat A
=(m,f_2,\hdots,f_r)\widehat A$. Therefore $P$ is a minimal prime
of $(f_1,\hdots,f_r)\widehat A$.

The extension
$$K[[Z_1]]\longrightarrow
K[[Z_1]]\,[Z_2]/(f_1,\hdots,f_r)$$ is integral and finite. Thus the
ring $K[[Z_1]]\,[Z_2]/(f_1,\hdots,f_r)$ is complete. This implies
$\widehat A=K[[x_1,z_2,\hdots,z_n]]$ is the
$(f_1,\hdots,f_r)$-adic (and the $(Z_2)$-adic) completion of
$K[[Z_1]]\,[Z_2]_{(Z)}$ and that $P$ is extended from
$K[[Z_1]]\,[Z_2]_{(Z)}$. This completes the proof of Proposition
2.2. \qed
\enddemo

The following theorem is the technical heart of the paper.

\proclaim{2.3 Theorem } Let $K$ be a field and let $y$ and
$X=\{x_1,\hdots,x_n\} $ be variables over  $K$. Assume that $V$ is
a discrete valuation domain with completion $\widehat V = K[[y]]$
and that $K[y]\subseteq V\subseteq K[[y]]$.  Also assume that the
field $K((y))=K[[y]]\,[1/y]$ has {\it uncountable} transcendence
degree over the quotient field ${\Cal Q}(V)$ of $V$. Set $R_0:=V[[X]]$
and $R = \widehat R_0=K[[y,X]]$. Let $P\in\Spec R$ be such that:
\roster \item"(i)" $P\subseteq (X )R$  (so $y\notin P$), and
\item"(ii)" $\dim(R/P) > 2.$
\endroster
Then there is a prime ideal $Q \in \Spec  R$ such that \roster
\item  $P\subset Q \subset XR, $
\item $\dim (R/Q) = 2$, and
\item  $P\cap R_0 = Q\cap R_0.$
\endroster
In particular, $P\cap K[[X]]= Q\cap K[[X]]$.
\endproclaim
\demo{Proof} Assume that $P$ has height $r$. Since $\dim(R/P) >
2$, we have  $0 \le r < n-1$. If $r > 0$, then there exist  a
transformation $x_1 \mapsto z_1, \hdots,x_n \mapsto z_n$ and
elements $f_1,\hdots,f_r\in P$, by Proposition 2.2, so that the
variable $y$ is fixed,  and

\medskip

 $f_1\in
K[[y,z_1,\hdots,z_{n-r}]]\,[z_{n-r+1},\hdots,z_n]\,$\text{ is monic
in
 $z_n$},

$f_2\in K[[y,z_1,\hdots,z_{n-r}]]\,[z_{n-r+1},\hdots,z_{n-1}]
$\text{ is monic in $z_{n-1}$ etc},

$\phantom{xx}$\vdots

 $f_r\in K[[y,z_1,\hdots,z_{n-r}]]\,[z_{n-r+1}]\, $\text{
is monic in $z_{n-r+1}$},

\medskip

 \noindent and the assertions of
Proposition 2.2  are satisfied. In particular, $P$ is a minimal
prime of $(f_1,\hdots,f_r)R$. Let $Z_1=\{z_1,\hdots,z_{n-r}\}$ and
$Z_2= \{z_{n-r+1},\hdots,z_{n-1},z_n\}.$ By Proposition 2.2, if $D
:= K[[y, Z_1]]\,[Z_2]_{(Z)}$ and $P_1 := P \cap D$, then $P_1R = P$.

The following diagram shows these rings and ideals.
 \vskip 18 pt

\setbox4=\vbox{\hbox{%
   \rlap{\kern2.5in\lower0in\hbox to .3in{\hss $R=K[[y,X]] = K[[y, Z_1, Z_2]]$ \hss}}%
   \rlap{\kern1.5in\lower0.8in\hbox to .3in{\hss $(X)R$ \hss}}%
  \rlap{\kern3.5in\lower1.1in\hbox to .3in{\hss $D =K[[y,Z_1]]\,[Z_2]_{(Z)}$
\hss}}%
    \rlap{\kern1.5in\lower1.55in\hbox to .3in{\hss $P=P_1R$ \hss}}%
   \rlap{\kern2.5in\lower2.3in\hbox to .3in{\hss $P_1 =P\cap D$ \hss}}%
   \rlap{\special{pa 2500 80} \special{pa 1500 650} \special{fp}}
   \rlap{\special{pa 2500 80} \special{pa 3500 950} \special{fp}}
       \rlap{\special{pa 1500  880} \special{pa 1500 1350} \special{fp}}
       \rlap{\special{pa 1500  1630} \special{pa 2500 2150} \special{fp}}
           \rlap{\special{pa 3500  1180} \special{pa 2500 2150}
\special{fp}}
    }}

\box4 \vskip 12 pt
\medskip

\noindent Note that $f_1,\ldots,f_r\in P_1$. Let $g_1, \ldots, g_s
\in P_1$ be other generators   such that

\noindent $ P_1= (f_1,\hdots,f_r,g_1,\hdots,g_s)D$. Then  $P= P_1R
= (f_1,\hdots,f_r,g_1,\hdots,g_s)R. $ For each $(i):=(i_1, \hdots,
i_n)\in \N^n$ and $j$, $k$ with $1\le j\le r$, $1\le k\le s$, let
$a_{j,(i)},\, b_{k,(i)}$ denote the coefficients in $K[[y]]$ of
the $ f_j,\, g_k$, so that
$$ f_j = \sum_{(i)\in {\Bbb N}^n} a_{j,(i)} z_1^{i_1}\hdots
z_n^{i_n}\, , \qquad g_k = \sum_{(i)\in {\Bbb N}^n} b_{k,(i)}
z_1^{i_1}\hdots z_n^{i_n}\in K[[y]]\,[[Z]].$$

Define
$$
\Delta := \cases \{a_{j,(i)}, b_{k,(i)}\}\subseteq K[[y]],
& \text{ for } r > 0 \\
\emptyset, & \text{ for } r = 0. \endcases
$$

\noindent A key observation  here is that in either case the set
$\Delta$ is countable.

To continue the proof, we consider $ S: = {\Cal Q}(V(\Delta))\cap
K[[y]]$, a discrete valuation domain, and its  field of quotients
$L:=\Cal Q(V(\Delta))$. Since $\Delta$ is a countable set, the
field $K((y))$ is (still) of uncountable transcendence degree over
$L$.  Let $\gamma_2,\hdots,\gamma_{n-r}$ be elements of $K[[y]]$
that are algebraically independent over $L$. We define

\noindent $T:=L(\gamma_2,\hdots,\gamma_{n-r})\cap K[[y]]$  and
$E:=\Cal Q(T)=L(\gamma_2,\hdots,\gamma_{n-r})$.

The diagram below shows the prime ideals $P$ and $P_1$ and the
containments between the relevant rings.

\vskip 18 pt

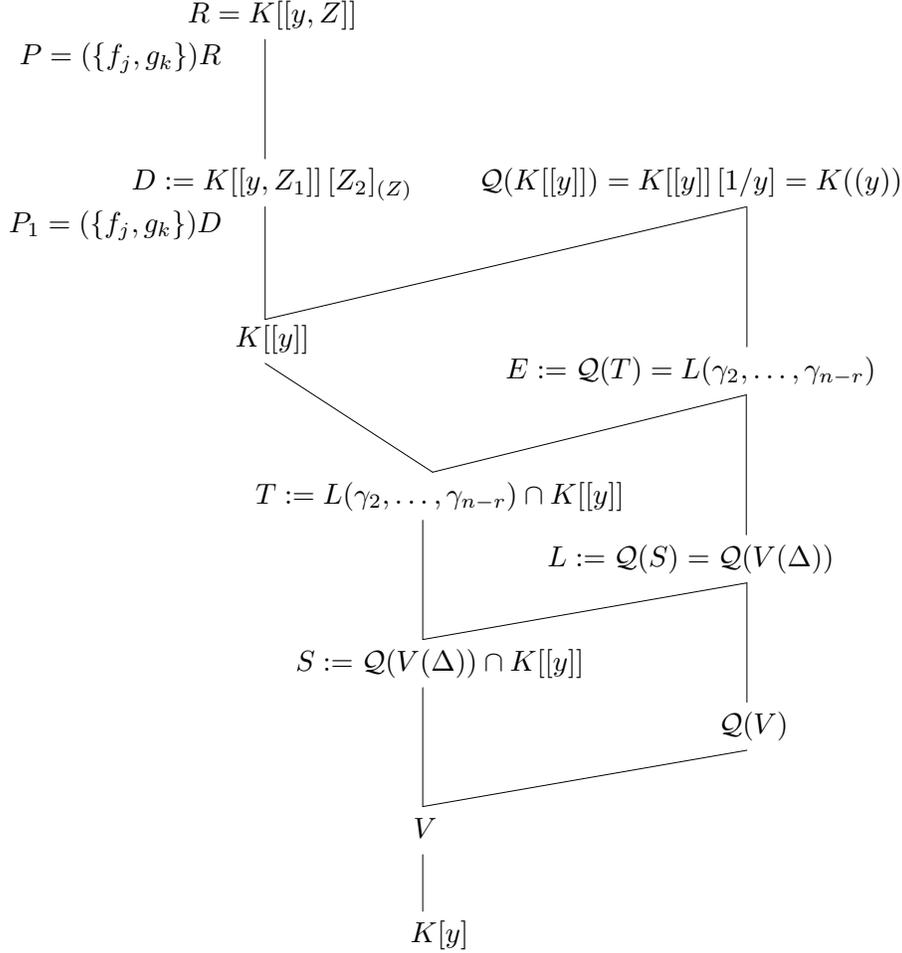
\setbox4=\vbox{\hbox{%
   \rlap{\kern1.5in\lower0in\hbox to .3in{\hss $R=K[[y,Z]]$ \hss}}%
    \rlap{\kern.75in\lower.2in\hbox to .3in{\hss $P=(\{f_j,g_k\})R$\hss}}%
     \rlap{\kern1.5in\lower0.8in\hbox to .3in{\hss
$D:=K[[y,Z_1]]\,[Z_2]_{(Z)}$ \hss}}%
   \rlap{\kern.75in\lower1.0in\hbox to .3in{\hss $P_1=(\{f_j,g_k\})D$
\hss}}%
     \rlap{\kern3.5in\lower.8in\hbox to .3in{\hss $\Cal
Q(K[[y]])=K[[y]]\,[1/y]=K((y))$ \hss}}%
    \rlap{\kern1.5in\lower1.55in\hbox to .3in{\hss $K[[y]]$ \hss}}%
   \rlap{\kern3.5in\lower1.70in\hbox to .3in{\hss $E:=\Cal
Q(T)=L(\gamma_2,\hdots,\gamma_{n-r})$ \hss}}%
    \rlap{\kern2.3in\lower2.3in\hbox to .3in{\hss
$T:=L(\gamma_2,\hdots,\gamma_{n-r})\cap K[[y]]$ \hss}}%
    \rlap{\kern2.3in\lower3.1in\hbox to .3in{\hss $S:=\Cal Q(V(\Delta))\cap
K[[y]]$ \hss}}%
      \rlap{\kern3.5in\lower2.6in\hbox to .3in{\hss $L:=\Cal Q(S)=\Cal
Q(V(\Delta))$ \hss}}%
   \rlap{\kern2.23in\lower3.9in\hbox to .3in{\hss $V$ \hss}}%
  \rlap{\kern3.8in\lower3.4in\hbox to .3in{\hss $\Cal Q(V)$ \hss}}%
   \rlap{\kern2.3in\lower4.4in\hbox to .3in{\hss $K[y]$ \hss}}%
           \rlap{\special{pa 1500 80} \special{pa 1500 650} \special{fp}}
       \rlap{\special{pa 1500 880} \special{pa 1500 1420} \special{fp}}
      \rlap{\special{pa 3800 880} \special{pa 1500 1420} \special{fp}}
        \rlap{\special{pa 3800 1780} \special{pa 3800 2450} \special{fp}}
       \rlap{\special{pa 3800 1780} \special{pa 2300 2150} \special{fp}}
                          \rlap{\special{pa 1500 1630} \special{pa 2300
2150} \special{fp}}
                 \rlap{\special{pa 2300 2380} \special{pa 2300
2950}\special{fp}}
                 \rlap{\special{pa 2300 3180} \special{pa 2300
3750}\special{fp}}
                   \rlap{\special{pa 3850  880} \special{pa 3850
1550}\special{fp}}
           \rlap{\special{pa 3850  2680} \special{pa 3850
3250}\special{fp}}
              \rlap{\special{pa 3850  2680} \special{pa 2300
2950}\special{fp}}
          \rlap{\special{pa 3850  3480} \special{pa 2300
3750}\special{fp}}
              \rlap{\special{pa 2300  3980} \special{pa 2300
4250}\special{fp}}
               }}

\box4 \vskip 12 pt

Let $P_2 := P \cap S[[Z_1]]\,[Z_2]_{(Z)}$. Since $f_1, \ldots, f_r,
g_1, \ldots, g_s  \in S[[Z_1]]\,[Z_2]_{(Z)}$, we have $P_2R = P$.
 Since $P\subseteq (x_1,\hdots,x_n)R = (Z)R$,
there is a prime ideal $\widetilde P$ in $L[[Z]]$ that is minimal
over $P_2L[[Z]]$. Since $L[[Z]]$ is flat over $S[[Z]]$,
$\widetilde P \cap S[[Z]] = P_2S[[Z]]$. Note that $L[[X]]=L[[Z]]$
is the $(f_1,\hdots,f_r)$-adic (and the $(Z_2)$-adic) completion
of $L[[Z_1]][Z_2]_{(Z)}$. In particular,
$$
L[[Z_1]]\,[Z_2]/(f_1,\hdots,f_r) =
 L[[Z_1]]\,[[Z_2]]/(f_1,\hdots,f_r)
$$
and this also holds with the field $L$ replaced by its extension
field $E$.

Since $L[[Z]]/\widetilde P$ is a homomorphic image of
$L[[Z]]/(f_1,\hdots,f_r)$,  it follows that $L[[Z]]/\widetilde P$
is integral (and finite) over $L[[Z_1]]$. This yields the
commutative diagram:

$$
\CD E[[Z_1]] &\longrightarrow &E[[Z_1]]\,[[Z_2]]/\widetilde P E[[Z]]
\\
\uparrow & &\uparrow \\
L[[Z_1]] &\longrightarrow &L[[Z_1]]\,[[Z_2]]/\widetilde P
\endCD\tag{2.3.0}
$$
with injective integral (finite) horizontal maps. Recall that
$E$ is the subfield of $K((y))$ obtained by adjoining $\gamma_2,
\ldots, \gamma_{n-r}$ to the field $L$. Thus the vertical maps
of (2.3.0) are faithfully flat.

Let $\q :=
(z_2-\gamma_2z_1,\hdots,z_{n-r}-\gamma_{n-r}z_1)E[[Z_1]]\in
\Spec(E[[Z_1]])$ and let $\widetilde W$ be a minimal prime of the
ideal $(\widetilde P, \q)E[[Z]]$.
 Since
$$
f_1,\hdots,f_r,z_2-\gamma_2z_1,\hdots,z_{n-r}-\gamma_{n-r}z_1
$$
is a regular sequence in $T[[Z]]$ the prime ideal $W:= \widetilde
W\cap T[[Z]]$ has height $n-1$. Let $\widetilde Q$ be a minimal
prime of $\widetilde W K((y))[[Z]]$ and let $Q:=\widetilde Q\cap
R$. Then  $W = Q  \cap T[[Z]]$, $P \subset Q  \subset ZR = XR$,  and
pictorially we have:

\vskip 18 pt

\setbox4=\vbox{\hbox{%
 \rlap{\kern2.1in\lower-.4in\hbox to .3in{\hss $(\widetilde
W)\subseteq\widetilde Q\subset K((y))[[Z]]$ \hss}}%
  \rlap{\kern1.5in\lower.1in\hbox to .3in{\hss $R:=K[[y,Z]]$ \hss}}%
  \rlap{\kern.9in\lower.5in\hbox to .3in{\hss $P=(\{f_j,g_k\})R \subseteq
Q\subset R$ \hss}}%
  \rlap{\kern1.5in\lower.9in\hbox to .3in{\hss $D:=K[[y,Z_1]]\,[Z_2]_{(Z)}$
\hss}}%
  \rlap{\kern.9in\lower1.1in\hbox to .3in{\hss $P_1=(\{f_j,g_k\})D\subset
D$\hss}}%
    \rlap{\kern3.2in\lower0in\hbox to .3in{\hss $ (\widetilde
P,\q)\subseteq\widetilde W\subset \phantom{x}E[[Z]]$ \hss}}%
    \rlap{\kern3.35in\lower.8in\hbox to .3in{\hss $(P_2)\subseteq\widetilde
P\subset L[[Z]]$ \hss}}%
\rlap{\kern1.2in\lower2.0in\hbox to .3in{\hss$\,P_2=(\{f_j,g_k\})
\subset
S[[Z_1]]\,[Z_2]_{(Z)}$ \hss}}%
     \rlap{\kern2.35in\lower.6in\hbox to .3in{\hss $\,W\subset
\phantom{x}T[[Z]]$ \hss}}%
       \rlap{\kern2.5in\lower1.3in\hbox to .3in{\hss
$\phantom{xx}S[[Z]]$\hss}}%
  \rlap{\kern3.6in\lower1.4in\hbox to .3in{\hss $L[[Z_1]]\,[Z_2]_{(Z)}$
\hss}}%
   \rlap{\kern4.2in\lower1.1in\hbox to .3in{\hss $\q\subset E[[Z_1]]$
\hss}}%
    \rlap{\kern4.3in\lower1.8in\hbox to .3in{\hss $L[[Z_1]]$ \hss}}%
          \rlap{\special{pa 2650 -320} \special{pa 1800 -50} \special{fp}}
         \rlap{\special{pa 1800 180} \special{pa 1800 750} \special{fp}}
              \rlap{\special{pa 1800 180} \special{pa 2650 450} \special{fp}}
          \rlap{\special{pa 2700 -320} \special{pa 3800 -150}\special{fp}}
          \rlap{\special{pa 1800 980} \special{pa 1800 1850} \special{fp}}
       \rlap{\special{pa 3750 80} \special{pa 3750 650} \special{fp}}
            \rlap{\special{pa 3800 880} \special{pa 3800
1250}\special{fp}}
                \rlap{\special{pa 3800 80} \special{pa 2650
450}\special{fp}}
        \rlap{\special{pa 2650 680} \special{pa 2650 1150} \special{fp}}
         \rlap{\special{pa 2650 1380} \special{pa 1800 1850} \special{fp}}
             \rlap{\special{pa 2700 1150} \special{pa 3800
880}\special{fp}}
\rlap{\special{pa 3750 1480} \special{pa 1800 1850} \special{fp}}
\rlap{\special{pa 3800 1480} \special{pa 4450 1680}\special{fp}}
\rlap{\special{pa 3800 80} \special{pa 4450 980}\special{fp}}
\rlap{\special{pa 4450 1150} \special{pa 4450 1680}\special{fp}}
               }}

\box4 \vskip 18 pt

Notice that $\q $  is a prime ideal of height $n-r-1$.  Also,
since $K((y))[[Z]]$ is flat over $K[[y, Z]] = R$, we have $\hgt Q
= n-1$ and $\dim(R/Q) = 2$.  We clearly have $P_2 \subseteq W \cap
S[[Z_1]]\,[Z_2]_{(Z)}$.

\medskip
\noindent
{\bf 2.3.1 Claim. }  $\q \cap L[[Z_1]] = (0)$.
\medskip

To show this we argue as in \cite{7}: Suppose that
$$ h = \sum_{m\in \Bbb N} H_m \, \in  \, \q\cap L[[z_1,\hdots,z_{n-r}]],$$
where $H_m \in L[z_1, \hdots , z_{n-r}]$  is a homogeneous
polynomial  of degree $m$:
$$
H_m = \sum_{|(i)|=m} c_{(i)} z_1^{i_1}\hdots z_{n-r}^{i_{n-r}},
$$
where $(i) := (i_1,\hdots,i_{n-r})\in \N^{n-r}$, $|(i)| := i_1 +
\cdots + i_{n-r}$  and $c_{(i)}\in L$. Consider the $E$-algebra
homomorphism $\pi : E[[Z_1]] \to E[[z_1]]$ defined by $\pi(z_1) =
z_1$ and $\pi(z_i) = \gamma_iz_1$ for $2 \le i \le n-r$. Then
$\ker \pi = \q$,  and  for each  $m\in \Bbb N$:
$$ \pi(H_m) = \pi(\sum_{|(i)|=m} c_{(i)} z_1^{i_1}\hdots z_{n-r}^{i_{n-r}}) =
\sum_{|(i)|=m} c_{(i)} \gamma_2^{i_2}\hdots \gamma_{n-r}^{i_{n-r}}z_1^m
$$
and
$$
\pi(h) = \sum_{m \in \N}\pi(H_m) =
\sum_{m \in \N}\sum_{|(i)|=m} c_{(i)} \gamma_2^{i_2}\hdots \gamma_{n-r}^{i_{n-r}}z_1^m.
$$
Since $h \in \q$, $\pi(h) = 0$. Since $\pi(h)$ is a power series in $E[[z_1]]$, 
each of its coefficients is zero, that is, for each $m \in \N$, 
$$
\sum_{|(i)|=m} c_{(i)} \gamma_2^{i_2}\hdots \gamma_{n-r}^{i_{n-r}} = 0.
$$
Since the $\gamma_i$ are algebraically independent over $L$,  each $c_{(i)} = 0$.
Therefore  $h=0$, and so $\q \cap L[[Z_1]] = (0)$. This proves Claim 2.3.1.

Using the commutativity of the displayed diagram (2.3.0) and  that
the horizonal maps of this diagram  are integral extensions, we
deduce that $(\widetilde W \cap E[[Z_1]] )= \q$, and $\q \cap
L[[Z_1]] = (0)$ implies $\widetilde W\cap L[[Z_1]] =(0)$. We
conclude that $Q \cap S[[Z]] = P \cap S[[Z]]$ and therefore $Q
\cap R_0 = P \cap R_0$. \qed
\enddemo

We record the following corollary.

 \proclaim{2.4 Corollary } Let $K$ be a field and let $R=K[[y,X]]$,
where 
$X=\{x_1,\hdots,x_n\} $ and $y$  are independent variables over  $K$.
Assume  $P\in \Spec R$ is such  that: \roster \item"(i)"
$P\subseteq (x_1,\hdots,x_n)R$ and \item"(ii)" $\dim(R/P) > 2.$
\endroster
Then there is a prime ideal $Q \in \Spec R$ so that \roster \item
$P\subset Q\subset (x_1,\hdots,x_n)R$, \item  $\dim (R/Q) = 2$,
and \item  $P\cap K[y]_{(y)}[[X]] = Q\cap K[y]_{(y)}[[X]]$.
\endroster
In particular, $P\cap K[[x_1,\hdots,x_n]]= Q\cap
K[[x_1,\hdots,x_n]]$.
\endproclaim
\demo{Proof}
With notation as in Theorem 2.3, let $V
= K[y]_{(y)}$.
\enddemo

\subheading{3. Weierstrass implications for the ring $B =
K[[X]]\,[Y]_{(X,Y)}$ }

As before $K$ denotes a field, $n$ and $m$ are positive integers,
and $X=\{x_1,\hdots,x_n\}$ and $Y=\{y_1,\hdots,y_m\}$ denote sets
of variables over $K$. Let  $B := K[[X ]]\, [Y]_{(X,Y)} =
K[[x_1,\hdots,x_n]]\,[y_1,\hdots,y_m]_{(x_1,\hdots,x_n,y_1,\hdots,y_m)}.$
The  completion of $B$ is  $\widehat B = K[[X,Y]]$.

 \proclaim{3.1
Theorem} With the notation as above, every ideal $Q$ of $\widehat
B = K[[X,Y]]$ maximal with the property that $Q\cap B =(0)$ is a
prime ideal of height $n+m-2$.
\endproclaim

\demo{Proof} Suppose first that $Q$ is such an ideal. Then clearly
$Q$ is prime. Matsumura shows in \cite{7, Theorem 3} that the
dimension of the generic formal fiber of $B$ is at most $n+m-2$.
Therefore $\hgt Q \le n+m-2$.

Now suppose  $P\in\Spec \widehat B$ is an arbitrary prime ideal of
height $r<n+m-2$ with $P\cap B=(0)$. We construct a prime
$Q\in\Spec \widehat B$ with $P\subset Q$, $Q\cap B= (0)$, and
$\hgt Q = n+m-2$. This will show that all prime ideals maximal in
the generic fiber have height $n+m-2$.

For the construction of $Q$ we consider first the case where $P
\not\subseteq X\widehat B$. Then  there exists a prime element
$f\in P$ that contains a term $\theta := y_1^{i_1}\cdots
y_m^{i_m}$, where the $i_j$'s are nonnegative integers and at
least one of the $i_j$ is positive. Notice that $m \ge 2$ for
otherwise with $y = y_1$ we have $f \in P$ contains a term $y^i$.
By Weierstrass it follows that $f = g\epsilon$, where $g \in
K[[X]]\,[y]$ is a nonzero  monic polynomial in $y$ and $\epsilon$ is
a unit of $\widehat B$. But $g \in P$ and $g \in B$ implies $P
\cap B \ne (0)$, a contradiction to our assumption that $P \cap B
= (0)$.

For convenience we now assume that the last exponent $i_m$
appearing in $\theta$ above is positive. We apply a change of
variables: $y_m\to t_m:=y_m$ and, for $1\le \ell< m$, let $
y_\ell\to t_\ell:=y_\ell+{t_m}^{e_\ell}$, where the $e_\ell$ are
chosen so that $f$, expressed in the variables $t_1,\hdots,t_m$,
contains a term $t_m^q$, for some positive integer $q$. This
change of variables induces an automorphism of $B$. By Weierstrass
$ f = g_1h,$ where $h$ is a unit in $\widehat B$ and $g_1\in
K[[X,t_1,\hdots,t_{m-1}]]\,[t_m]$ is monic in $t_m$. Set $P_1= P\cap
K[[X,t_1,\hdots,t_{m-1}]]$. If $P_1\subseteq XK[[X, t_1, \ldots,
t_{m-1}]]$, we stop the procedure and take $s=m-1$ in what
follows. If $P_1\not\subseteq XK[[X, t_1, \ldots, t_{m-1}]]$, then
there exists a prime element $\tilde f\in P_1$ that contains a
term ${t_1}^{j_1}\cdots {t_{m-1}}^{j_{m-1}}$, where the $j_k$'s
are nonnegative integers and at least one of the $j_k$ is
positive. We then  repeat the procedure using the prime ideal
$P_1$. That is, we replace $ t_1,\dots t_{m-1}$ with a change of
variables so that a prime element of $P_1$ contains a term monic
in some one of the new variables. After a suitable finite
iteration of changes of variables, we obtain an automorphism of
$\widehat B$ that restricts to an automorphism of $B$ and  maps
$y_1,\hdots,y_m\mapsto z_1,\hdots,z_m$. Moreover, there exist a
positive integer $s \le m-1$ and elements $g_1, \ldots g_{m-s} \in
P$ such that

$$\alignat2
g_1&\in K[[X, z_1,\hdots,z_{m-1}]]\,[z_m]\quad&&\text{is monic in
$z_m$}\\
g_2&\in K[[X, z_1,\hdots,z_{m-2}]]\,[z_{m-1}] \quad&&\text{is monic
in $z_{m-1}$, etc} \\
 &\vdots && \\
g_{m-s}&\in K[[X, z_1,\hdots,z_s]]\,[z_{s+1}] \quad&&\text{is monic
in $z_{s+1}$,}
\endalignat
$$
and such that, for $R_s:= K[[X, z_1, \ldots,z_s]]$ and $P_s := P
\cap R_s$, we have $P_s\subseteq XR_s$.

As in the proof of Proposition 2.2 we replace the regular sequence
$g_1,\hdots.g_{m-s}$ by a regular sequence $f_1,\hdots,f_{m-s}$ so
that:
$$\alignat2
f_1&\in R_s[z_{s+1},\hdots,z_m]\quad&&\text{is monic in
$z_m$}\\
f_2&\in R_s[z_{s+1},\hdots,z_{m-1}]
\quad&&\text{is monic in $z_{m-1}$, etc}\\
 &\vdots && \\
f_{m-s}&\in R_s[z_{s+1}] \quad&&\text{is monic in $z_{s+1}$.}
\endalignat
$$
and $(g_1,\hdots.g_{m-s})\widehat B = (f_1,\hdots,f_{m-s})\widehat
B$.

Let $G := K[[X, z_1, \ldots, z_s]]\,[z_{s+1}, \ldots,
z_m]=R_s[z_{s+1}, \ldots, z_m]$.  By
 Proposition 2.2, $P$ is extended from $G$. Let
 $\q := P \cap G$ and
extend $f_1,\hdots,f_{m-s}$ to a generating system of $\q$, say,
$\q = (f_1,\hdots,f_{m-s},h_1,\hdots,h_t)G$. For integers
$k,\,\ell$ with $1\le k\le m-s$ and $1\le \ell\le t$,  express the
$f_k$ and $h_\ell$ in $G$ as power series in $\widehat
B=K[[z_1]][[z_2, \ldots, z_m]]\,[[X]]$ with coefficients in
$K[[z_1]]$:
$$
f_k = \sum a_{k(i)(j)} z_2^{i_2}\hdots z_m^{i_m}x_1^{j_1}\hdots
x_n^{j_n} \quad \text{and}   \quad h_\ell = \sum b_{\ell(i)(j)}
z_2^{i_2}\hdots z_m^{i_m}x_1^{j_1}\hdots x_n^{j_n}\,,
$$
where  $a_{k(i)(j)}, b_{\ell(i)(j)}\in K[[z_1]]$, $(i) = (i_2,
\ldots, i_m)$  and $(j) = (j_1, \ldots, j_n)$. The set $\Delta =
\{a_{k(i)(j)}, b_{\ell(i)(j)}\}$ is countable. We define $V :
=K(z_1,\Delta)\cap K[[z_1]]$. Then $V$ is a discrete valuation
domain with completion $K[[z_1]]$ and $K((z_1))$ has  uncountable
transcendence degree over $\Cal Q(V)$. Let
$V_s:=V[[X,z_2,\hdots,z_s]]\subseteq R_s$. Notice that $R_s
=\widehat{V_s}$, the completion of $V_s$. Also $f_1,\dots,f_{
m-s}\in V_s[z_{s+1}, \ldots, z_m] \subseteq G$ and $(f_1,\dots,f_{
m-s})G\cap R_s=(0)$. Furthermore the extension
$$
V_s:=V[[X,z_2,\hdots,z_s]]\hookrightarrow
 V_s[z_{s+1},\hdots,z_m]/(f_1,\hdots,f_{m-s})
$$
is finite. Set $P_0 := P\cap V_s$. Then $P_0\subseteq XR_s \cap
V_s = XV_s$.

Consider the commutative diagram:
$$
\CD R_s:=K[[X,z_1,\hdots,z_s]] &\longrightarrow
&R_s[[z_{s+1},\hdots,z_m]]/(f_1,\hdots,f_{m-s})\\
\uparrow & &\uparrow \\
 V_s:=V[[X,z_2,\hdots,z_s]] &\longrightarrow
&V_s[z_{s+1},\hdots,z_m]/(f_1,\hdots,f_{m-s})\,.
\endCD \tag{3.1.1}
$$
The horizontal maps are injective and finite and the vertical maps
are completions.

The prime ideal $\bar {\q} :=
PR_s[[z_{s+1},\hdots,z_m]]/(f_1,\hdots,f_{m-s})$ lies over $P_s$
in $R_s$. By assumption $P_s\subseteq (X)R_s$ and by Theorem 2.3
there is a prime ideal $Q_s$ of $R_s$ such that  $P_s\subseteq Q_s
\subseteq (X)R_s$, \, $ Q_s \cap V_s = P_s\cap V_s=P_0$, and $\dim
(R_s/ Q_s) =2$. There is a prime ideal $\bar Q$ in
$R_s[[z_{s+1},\hdots,z_m]]/(f_1,\hdots,f_{m-s})$ lying over $ Q_s$
with $\bar \q\subseteq \bar Q$ by the ``going-up theorem''
\cite{8, Theorem 9.4}. Let $Q$ be the preimage in $\widehat B =
K[[X,z_1,\hdots,z_m]]$ of $\bar Q$. We show the rings and ideals
of Theorem 3.1 below.

\vskip 18 pt

\setbox4=\vbox{\hbox{%
   \rlap{\kern2.5in\lower0in\hbox to .3in{\hss $\widehat B=K[[X,Y]]=K[[X,z_1,\ldots, z_m]]=R_s[[z_{s+1},\ldots, z_m]]$ \hss}}%
    \rlap{\kern1.75in\lower.2in\hbox to .3in{\hss $(\q,Q_s)\widehat B\subseteq Q$\hss}}%
    \rlap{\kern1.5in\lower.4in\hbox to .3in{\hss $P\nsubseteq X\widehat B$\hss}}%
         \rlap{\kern2.5in\lower0.8in\hbox to .3in{\hss
$G:=R_s[z_{s+1},\ldots,z_m]$ \hss}}%
   \rlap{\kern2.0in\lower1.0in\hbox to .3in{\hss $\q:=P\cap G$
\hss}}%
 \rlap{\kern2.0in\lower1.2in\hbox to .3in{\hss $\q=(\{f_i,h_j\})G$
\hss}}%
    \rlap{\kern2.4in\lower1.55in\hbox to .3in{\hss $f_i\notin R_s:=K[[X,z_1,\ldots,z_s]]$ \hss}}%
     \rlap{\kern1.8in\lower1.75in\hbox to .3in{\hss $P_s\subseteq Q_s\subset R_s$ \hss}}%
      \rlap{\kern1.0in\lower1.95in\hbox to .3in{\hss $P_s:=P\cap R_s\subseteq XR_s$ \hss}}%
       \rlap{\kern3.5in\lower2.3in\hbox to .3in{\hss $\widehat V
=K[[z_1]]$ \hss}}%
    \rlap{\kern1.5in\lower2.3in\hbox to .3in{\hss
$V_s:=V[[X,z_2,\ldots,z_s]]$ \hss}}%
   \rlap{\kern1.050in\lower2.5in\hbox to .3in{\hss
$P_0:=P\cap V_s$ \hss}}%
    \rlap{\kern2.5in\lower3.1in\hbox to .3in{\hss $V:=K(z_1,\Delta)\cap
K[[z_1]]$ \hss}}%
           \rlap{\special{pa 2580 80} \special{pa 2580 650} \special{fp}}
            \rlap{\special{pa 2580 1630} \special{pa 3580 2150} \special{fp}}
                \rlap{\special{pa 2580 880} \special{pa 2580 1420} \special{fp}}
                          \rlap{\special{pa 2580 1630} \special{pa 1500
2150} \special{fp}}
                 \rlap{\special{pa 1500 2380} \special{pa 2580
2950}\special{fp}}
              \rlap{\special{pa 3580  2380} \special{pa 2580
2950}\special{fp}}
               }}

\box4

\vskip 18 pt

Then $Q$ has height $n+s-2 + m-s = n+m-2$. Moreover,  from diagram
(3.1.1),  it follows  that $Q$ and $P$ have the same contraction
to $V_s[z_{s+1},\hdots,z_m]$. This implies that $Q\cap B= (0)$ and
completes the proof in the case where $P \not \subseteq X\widehat
B$.

In the case where  $P\subseteq X\widehat B$, let $h_1, \ldots, h_t
\in \widehat B$ be a finite set of generators of $P$, and as
above, let $b_{\ell(i)(j)} \in K[[z_1]]$ be the coefficients of
the $h_\ell$'s. Consider the countable set $\Delta = \{
b_{\ell(i)(j)}\}$ and the valuation domain $V := K(z_1,\Delta)\cap
K[[z_1]]$. Set $P_0 := P \cap V[[X, z_2, \ldots, z_m]]$. By
Theorem 2.3,   there exists  a prime ideal $Q$ of $\widehat B =
K[[X, z_1, \ldots, z_m]]$ of height $n+m-2$  such that $P\subset
Q$ and  $Q \cap V[[X,z_2,\hdots,z_m]] =P\cap V[[X,z_2,\hdots,z_m]]
= P_0$. Therefore $Q \cap B = (0)$. This completes the proof of
Theorem 3.1. \qed
\enddemo

\subheading{4. Weierstrass implications for the ring $C
=K[Y]_{(Y)}[[X]]$}

As before $K$ denotes a field, $n$ and $m$ are positive integers,
and $X=\{x_1,\hdots,x_n\}$ and $Y=\{y_1,\hdots,y_m\}$ denote sets
of variables over $K$.  Consider the ring $C
=K[y_1,\hdots,y_m]_{(y_1,\hdots,y_m)}[[x_1,\hdots,x_n]]
=K[Y]_{(Y)}[[X]]$. Then the completion of $C$ is  $\widehat C =
K[[Y,X]]$.

 \proclaim{4.1 Theorem} With notation as above, let $Q\in\Spec  \widehat C $
be maximal with the property that $Q\cap C = (0)$. Then $\hgt Q =
n+m-2$.
\endproclaim

\demo{Proof} Let $B = K[[X]]\,[Y]_{(X,Y)} \subset C$. If $P \in
\Spec \widehat C = \Spec \widehat B$ and $P \cap C = (0)$, then $P
\cap B = (0)$, so $\hgt P \le n+m-2$ by Theorem~3.1.
 Consider a nonzero prime  $P\in\Spec \widehat
C$ with $P\cap C =(0)$ and $\hgt P = r < n+m-2$. If $P \subseteq
X\widehat C$ then Theorem 2.3 implies the existence of $Q \in
\Spec \widehat C$ with $\hgt Q = n+m-2$ such that $P \subset Q$
and $Q \cap C = (0)$.

Assume that $P$ is not contained in $X\widehat C$ and consider the
ideal $J := (P, X)\widehat C$. Since $C$ is complete in the
$XC$-adic topology, \cite{9, Lemma 2} implies that if $J$ is
primary for the maximal ideal of $\widehat C$, then $P$ is
extended from $C$. Since we are assuming $P \cap C = (0)$, $J$ is
not primary for the maximal ideal of $\widehat C$ and we have
$\hgt J =  n+s < n+m$, where $0 < s < m$. Let $W \in \Spec
\widehat C$ be a minimal prime of $J$ such that $\hgt W = n + s$.
Let $W_0 = W \cap K[[Y]]$. Then $W=(W_0,X)\widehat C$ and $W_0$ is
a prime ideal of $K[[Y]]$ with $\hgt W_0 = s$. By Proposition 2.2
applied to $K[[Y]]$ and the prime ideal $W_0 \in \Spec K[[Y]]$,
there exists a change of variables $Y\mapsto Z$ with $y_1\mapsto
z_1,\hdots,y_m\mapsto z_m$ and elements $f_1,\hdots,f_s\in W_0$ so
that with $Z_1 = \{z_1, \ldots, z_{m-s} \}$, we have
$$\alignat2
f_1&\in K[[Z_1]]\,[z_{m-s+1},\hdots,z_m]\quad&&\text{is monic in
$z_m$}\\
f_2&\in K[[Z_1]]\,[z_{m-s+1},\hdots,z_{m-1}] \quad&&\text{is monic
in $z_{m-1}$, etc}\\
 & \vdots && \\
f_s&\in K[[Z_1]]\,[z_{m-s+1}] \quad&&\text{is monic in
$z_{m-s+1}$.}
\endalignat
$$
Now $z_1,\hdots,z_{m-s},f_1,\hdots,f_s$ is a regular sequence in
$K[[Z]] = K[[Y]]$. Let $T = \{t_{m-s+1},\hdots,t_m \} $ be a set
of  additional variables and consider the map:
$$\varphi:K[[ Z_1, T]]\longrightarrow
K[[z_1,\hdots,z_m]]$$ defined by $z_i\mapsto z_i$ for all $1\le
i\le m-s$ and $t_{m-i+1}\mapsto f_i$ for all $1\le i\le s$. The
embedding $\varphi$ is finite (and free) and so is the extension
to power series rings in $X$:

$$\rho:K[[ Z_1, T]]\,[[X]]\longrightarrow
K[[z_1,\hdots,z_m]]\,[[X]] = \widehat C.
$$
Since $W \in \Spec \widehat C$ is  of height $n+s$, so is  its
contraction $\rho^{-1}(W) \in \Spec K[[Z_1, T, X]]$. Moreover
$\rho^{-1}(W)$ contains $(T, X)K[[Z_1, T, X]]$, a prime ideal of
height $n+s$. Therefore $\rho^{-1}(W) =(T, X)K[[Z_1, T, X]]$. By
construction, $P\subseteq W$ which yields that
$\rho^{-1}(P)\subseteq (T, X)K[[Z_1, T, X]]$.

To complete the proof we construct a suitable base ring related to
$C$. Consider the expressions for  the $f_i$'s as power series in
$z_2,\hdots,z_m$ with coefficients in $K[[z_1]]$:
$$ f_j = \sum a_{j(i)}z_2^{i_2}\hdots z_m^{i_m},$$
where $(i):=(i_2,\dots,i_m), 1\le j\le s, a_{j(i)}\in K[[z_1]]$.
Also consider a finite generating system $g_1,\hdots,g_q$ for $P$
and expressions for the $g_k$, where $1\le k\le q$, as power
series in $z_2, \ldots, z_m, x_1, \ldots, x_n$ with coefficients
in $K[[z_1]]$:
$$ g_k = \sum b_{k(i)(\ell)}z_2^{i_2}\hdots z_m^{i_m}x_1^{\ell_1}\hdots
x_n^{\ell_n},
$$
where $(i):=(i_2,\dots,i_m),\, (\ell):=(\ell_1,\dots,\ell_n),$ and
$ b_{k(i)(\ell)}\in K[[z_1]]$. We take the subset $\Delta =
\{a_{j(i)},\, b_{k(i)(\ell)} \}$ of $K[[z_1]]$ and consider the
discrete valuation domain:
$$ V :=K(z_1, \Delta )\cap K[[z_1]].$$
Since $V$ is countably generated over $K(z_1)$, the field $K((z_1))$ has
uncountable  transcendence degree over $\Cal Q(V) =
K(z_1,\Delta)$. Moreover, by construction the ideal $P$ is
extended from $V[[z_2,\hdots,z_m]]\,[[X]]$. Consider the embedding:
$$\psi:V[[z_2,\hdots,z_{m-s}, T ]]\longrightarrow
V[[z_2,\hdots,z_m]], $$ which is the restriction of $\varphi$
above, so that $z_i\mapsto z_i$ for all $2\le i\le m-s$ and
$t_{m-i+1}\mapsto f_i$ for all $i$ with $1\le i\le s$.

Let $\sigma$ be the extension of $\psi$ to the power series rings:

$$\sigma:V[[z_2,\hdots,z_{m-s}, T ]]\,[[X]]\longrightarrow
V[[z_2,\hdots,z_m]]\,[[X]]$$ with $\sigma(x_i) = x_i$ for all $i$ with 
$1\le i\le n$.

Notice that $\rho$ defined above is the completion $\widehat \sigma$ 
of the map $\sigma$, that is, 
the extension of $\sigma$ to the
completions. Consider the commutative diagram:
$$
\CD
K[[Z_1, T]]\,[[X]] @>{\widehat \sigma=\rho}>> K[[Z]]\,[[X]] = \widehat C \\
@AAA        @AAA  \\
V[[z_2,\hdots,z_{m-s}, T]]\,[[X]] @>{\sigma}>>
V[[z_2,\hdots,z_{m}]]\,[[X]]
\endCD \tag{4.1.0}
$$
where $\widehat \sigma=\rho$ is a finite map.

Recall that $\rho^{-1}(W) =(T, X)K[[Z_1, T, X]]$, and so $
\rho^{-1}(P)\subseteq (T, X)K[[Z_1, T, X]]$ by Diagram 4.1.0. By Theorem
2.3, there exists  a prime ideal $Q_0$ of the ring $K[[Z_1, T,
X]]$ such that $\rho^{-1}(P)\subseteq Q_0$, $\hgt Q_0  = n+m-2$,
and
$$
Q_0\cap V[[z_2,\hdots,z_{m-s}, T ]]\,[[X]] = \rho^{-1}(P)\cap
V[[z_2,\hdots,z_{m-s}, T ]]\,[[X]]. $$ By the ``going-up theorem''
\cite{8, Theorem 9.4}, there is a prime ideal $Q\in \Spec \widehat
C$ that lies over $Q_0$ and contains $P$. Moreover, $Q$ also has
height $n+m-2$. The commutativity of diagram (4.1.0) implies that

$$ P_1:= P \cap V[[z_2,\hdots,z_{m-s},T]]\,[[X]]\,\,
\subseteq \,\,Q_1:=Q \cap  V[[z_2,\hdots,z_{m-s},T]]\,[[X]].
$$

Consider the finite homomorphism:
$$\lambda: V[[z_2,\hdots,z_{m-s}]]\,[T]_{(Z_1,T)}[[X]]\longrightarrow
V[[z_2,\hdots,z_{m-s}]]\,[z_{m-s+1},\hdots,z_m]_{(Z)}[[X]]$$
(determined by $t_i\mapsto f_i$ for $1\le i\le m$) and the
commutative diagram:

$$\CD
V[[z_2,\hdots,z_{m-s}]]\,[[T]]\,[[X]] @>{\sigma}>> V[[z_2,\hdots,z_m]]\,[[X]] \\
@AAA     @AAA \\
V[[z_2,\hdots,z_{m-s}]]\,[T]_{(Z_1,T)}[[X]] @>{\lambda}>>
V[[z_2,\hdots,z_{m-s}]]\,[z_{m-s+1},\hdots,z_m]_{(Z)}[[X]].
\endCD
$$

\smallskip

Since $Q\cap V[[z_2,\hdots,z_{m-s},T]]\,[[X]]=P\cap
V[[z_2,\hdots,z_{m-s},T]]\,[[X]]$ and since $\lambda$ is a finite
map we conclude that
$$\multline \phantom{xxx}Q_1\cap V[[z_2,\hdots,z_{m-s}]]\,[z_{m-s+1},\hdots,z_m]_{(Z)}[[X]]\\
= P_1\cap
V[[z_2,\hdots,z_m]]\,[z_{m-s+1},\hdots,z_m]_{(Z)}[[X]].\phantom{xxx}\endmultline$$
Since $C\subseteq
V[[z_2,\hdots,z_{m-s}]]\,[z_{m-s+1},\hdots,z_m]_{(Z)}[[X]]$,  we
obtain that $Q \cap C =$ \linebreak $P \cap C = (0)$. This
completes the proof of Theorem 4.1.\qed
 \enddemo

\subheading{4.2 Remark} With $B$ and $C$ as in Sections 3 and 4,
we have
$$
B = K[[X]]\,[Y]_{(X, Y)} \hookrightarrow K[Y]_{(Y)}[[X]] = C \quad
\text{  and  } \quad \widehat B = K[[X, Y]] = \widehat C.
$$
Thus for $P \in \Spec K[[X, Y]]$, if $P \cap C = (0)$, then $P
\cap B = (0)$.  By Theorems 3.1 and 4.1, each prime of $K[[X, Y]]$
maximal in the generic formal fiber of $B$ or $C$ has height $n +
m - 2$. Therefore each $P \in \Spec K[[X, Y]]$ maximal with
respect to $P \cap C = (0)$ is also maximal with respect to $P
\cap B = (0)$. However, if $n+m \ge 3$, the generic fiber of $B
\hookrightarrow C$ is nonzero \cite{4}, so there exist  primes of
$K[[X, Y]]$ maximal in the generic formal fiber of $B$ that are
not in the generic formal fiber of $C$.

\subheading{5. Subrings of the power series ring $K[[z, t]]$}

In this section we establish properties of certain subrings of the
power series ring $K[[z, t]]$ that will be useful in considering
the generic formal fiber of localized polynomial rings over the
field $K$.

 \subheading{5.1. Notation} Let $K$ be a field and let
$z$ and $t$ be independent variables over $K$. Consider countably
many power series:
$$ \alpha_i(z) = \sum_{j=0}^{\infty} a_{ij} z^j\in K[[z]]$$
with coefficients $a_{ik}\in K$.  Let $s$ be a positive integer
and let  $\omega_1,\hdots,\omega_s\in K[[z, t]]$ be power series
in $z$ and $t$, say:
$$
\omega_i = \sum_{j=0}^{\infty} \beta_{ij}t^j, \quad\text{where}
\quad \beta_{ij}(z) = \sum_{k=0}^{\infty} b_{ijk}z^k\in K[[z]]
\quad \text{and}\quad b_{ijk} \in K,
$$
for each $i$ with $1 \le i \le s$. Consider the subfield $K(z,
\{\alpha_i\}, \{\beta_{ij}\} )$ of $K((z))$ and the discrete
rank-one valuation domain
$$ V := K(z, \{\alpha_i\}, \{\beta_{ij}\}) \cap K[[z]].$$
The completion of  $V$ is  $\widehat V=K[[z]]$. Assume that
$\omega_1,\hdots,\omega_r$ are algebraically independent over
$\Cal Q(V)(t)$ and that the elements
$\omega_{r+1},\hdots,\omega_s$ are algebraic over the field $\Cal
Q(V)(t, \{\omega_i\}_{i=1}^r)$. Notice that the set $\{\alpha_i\}
\cup \{\beta_{ij}\}$ is countable, and that also the set of
coefficients of the $\alpha_i$ and $\beta_{ij}$
$$\Delta :=\{a_{ij}\} \cup \{b_{ijk}\}$$
is a countable subset of the field $K$.  Let $K_0$ denote the
prime subfield of $K$ and let $F$ denote the algebraic closure in
$K$ of the field $K_0(\Delta)$. The field $F$ is countable and the
power series  $\alpha_i(z)$ and $\beta_{ij}(z)$ are in $F[[z]]$.
Consider the subfield $F(z, \{\alpha_i\}, \{\beta_{ij}\})$ of
$F((z))$ and the discrete rank-one valuation domain
$$ V_0 := F(z, \{\alpha_i\}, \{\beta_{ij}\})\cap  F[[z]].$$
The completion of  $V_0$ is  $\widehat V_0 = F[[z]]$. Since $\Cal
Q(V_0)(t) \subseteq \Cal Q(V)(t)$, the elements
$\omega_1,\hdots,\omega_r$ are algebraically independent over the
field $\Cal Q(V_0)(t)$.

Consider the subfield $E_0 := \Cal Q(V_0)(t, \omega_1, \ldots,
\omega_r) $ of $\Cal Q(V_0[[t]])$ and the subfield $E := \Cal
Q(V)(t, \omega_1, \ldots, \omega_r)$ of $\Cal Q(V[[t]])$. A result
of Valabrega \cite{11} implies that the integral domains:
$$
D_0 := E_0\cap V_0[[t]] \qquad \text{and} \qquad D := E\cap V[[t]]
\tag{5.1.1}
$$
are two-dimensional regular local rings with completions
$\widehat D_0 = F[[z, t]]$ and $\widehat D = K[[z, t]]$,
respectively. Moreover, $\Cal Q(D_0) = E_0$  is a countable field.

\proclaim{5.2 Proposition} Let $D_0$ be as defined in (5.1). Then
there exists  a power series $\gamma\in zF[[z]]$ such that the
prime ideal $(t-\gamma)F[[z, t]] \cap D_0 = (0)$, i.e., $(t -
\gamma)F[[z, t]]$ is in the generic formal fiber of $D_0$.
\endproclaim \demo{Proof} Since $D_0$ is
countable there are only countably many prime ideals in $D_0$ and
since $D_0$ is Noetherian there are only countably many prime
ideals in $\widehat D_0 =F[[z, t]]$ that lie over a nonzero prime
of $D_0$. There are uncountably many primes in $F[[z, t]]$, which
are generated by elements of the form $t -\sigma$ for some
$\sigma\in zF[[z]]$. Thus there must exist an element  $\gamma \in
zF[[z]]$ with $(t - \gamma)F[[z, t]] \cap D_0 =(0)$. \qed
\enddemo

For  $\omega_i = \omega_i(t) = \sum_{j=0}^\infty\beta_{ij}t^j$ as
in (5.1) and $\gamma$ an element of $zK[[z]]$,  let
$\omega_i(\gamma)$ denote the following power series in $K[[z]]$:
$$\omega_i(\gamma) := \sum_{j=0}^{\infty} \beta_{ij} \gamma^j\in K[[z]].$$

\proclaim{5.3 Proposition} Let $D$ be as defined in (5.1.1). For
an element $\gamma\in zK[[z]]$ the following conditions are
equivalent: \roster \item"(i)" $(t-\gamma)K[[z, t]]\cap D = (0)$.
\item"(ii)" The elements
$\gamma,\omega_1(\gamma),\hdots,\omega_r(\gamma)$ are
algebraically independent over $\Cal Q(V)$. \endroster
\endproclaim
\demo{Proof} $($i$)$ $\Rightarrow$ $($ii$)$: Assume by way of
contradiction that the set
$\{\gamma,\omega_1(\gamma),\hdots,\omega_r(\gamma)\}$ is
algebraically dependent over $\Cal Q(V)$ and let $d_{(k)} \in V$
be finitely many elements such that
$$
\sum_{(k)} d_{(k)} \omega_1(\gamma)^{k_{1}}\hdots
\omega_r(\gamma)^{k_{r}}\gamma^{k_{{r+1}}} = 0$$  is a nontrivial
equation of algebraic dependence for
$\gamma,\omega_1(\gamma),\hdots,\omega_r(\gamma)$, where each
$(k)=(k_{1}\hdots ,k_{r},k_{{r+1}})$ is an $(r+1)$-tuple of
nonnegative integers. It follows that
$$ \sum_{(k)} d_{(k)} \omega_1^{k_1}\hdots
\omega_r^{k_r}t^{k_{r+1}}\in (t-\gamma)K[[z, t]] \cap D = (0).$$
Since
$\omega_1,\hdots,\omega_r$ are algebraically independent over
$\Cal Q(V)(t)$,  we have $d_{(k)}=0$ for all $(k)$, a contradiction.  
This completes the proof that $(i)
\Rightarrow (ii)$.

$($ii$)$ $\Rightarrow$ $($i$)$: If $(t-\gamma)K[[z, t]]\cap D\ne
(0)$,  then there exists  a nonzero  element $$\tau =\sum_{(k)}
d_{(k)} \omega_1^{k_1}\hdots \omega_r^{k_r}t^{k_{r+1}}\in
(t-\gamma)K[[z, t]]\cap V[t,\omega_1,\hdots,\omega_r].$$ But this
implies that
$$\tau(\gamma) =\sum_{(k)} d_{(k)} \omega_1(\gamma)^{k_1}\hdots
\omega_r(\gamma)^{k_r}\gamma^{k_{r+1}} = 0.$$ Since
$\gamma,\omega_1(\gamma),\hdots,\omega_r(\gamma)$ are
algebraically independent over $\Cal Q(V)$, it follows that  all
the coefficients $d_{(k)}=0$, a contradiction to the  assumption
that $\tau$ is nonzero. \qed
\enddemo

\medskip

Let $\gamma\in zF[[z]]$ be as in Proposition 5.2 with
$(t-\gamma)F[[z, t]] \cap D_0 = (0)$. Then:

\proclaim{5.4 Proposition} With notation as above,  we also have
$(t-\gamma)K[[z, t]]\cap D = (0)$, that is, $(t-\gamma)K[[z, t]]$
is in the generic formal fiber of $D$.
\endproclaim

\demo{Proof} Let $\{t_i\}_{i\in I}$ be a transcendence basis of
$K$ over $F$ and let  $L := F(\{t_i\}_{i\in I})$. Then $K$ is
algebraic over $L$. Let  $\{\alpha_i\}, \{\beta_{ij}\} \subset
F[[z]]$  be as in (5.1) and define
$$ V_1 = L(z, \{\alpha_i\}, \{\beta_{ij}\})\cap L[[z]]\quad \text{and} \quad  D_1 =
\Cal Q(V_1)(t, \omega_1,\hdots,\omega_r)\cap L[[z, t]].$$ Then
$V_1$ is a discrete rank-one valuation domain with completion
$L[[z]]$ and $D_1$ is a two-dimensional  regular local domain with
completion $\widehat D_1 = L[[z, t]]$. Note that $\Cal Q(V)$ and
$\Cal Q(D)$ are algebraic over $\Cal Q(V_1)$ and $\Cal Q(D_1)$,
respectively. Since $(t-\gamma)K[[z, t]]\cap L[[z, t]] =
(t-\gamma)L[[z, t]]$, it suffices to prove that $(t-\gamma)L[[z,
t]] \cap D_1 =(0)$. By Proposition 5.3, it suffices to show that
$\gamma,\omega_1(\gamma),\hdots,\omega_r(\gamma)$ are
algebraically independent over $\Cal Q(V_1)$. The  commutative
diagram

$$
\CD
F[[z]] @>{\{t_i\} {\text{algebraically ind.}} }>> L[[z]]\\
@AAA             @AAA \\
\Cal Q(V_0) @>{\text{transcendence basis } \{t_i\}}>>  \Cal Q(V_1)
\endCD
$$
implies  that the set
$\{\gamma,\omega_1(\gamma),\hdots,\omega_r(\gamma)\} \cup \{t_i\}$
is algebraically independent over $\Cal Q(V_0)$. Therefore
$\{\gamma,\omega_1(\gamma),\hdots,\omega_r(\gamma)\}$ is
algebraically independent over $\Cal Q(V_1)$, which completes the
proof of Proposition 5.4. \qed
\enddemo

\subheading{5.5 Remark}  We remark that with
$\omega_{r+1},\hdots,\omega_s$ algebraic over $\Cal
Q(V)(\omega_1,\hdots,\omega_r)$ as in (5.1), if we define
$$\widetilde D := \Cal Q(V)(t, \omega_1,\hdots,\omega_s)\cap V[[t]],$$
then again by Valabrega \cite{11},  $\widetilde D$ is a
two-dimensional regular local domain  with completion $K[[z, t]]$.
Moreover, $\Cal Q(\widetilde D)$ is algebraic over $\Cal Q(D)$ and
$(t-\gamma)K[[z, t]] \cap D = (0)$ implies that $(t-\gamma)K[[z,
t]] \cap \widetilde D =(0)$.

\subheading{6. Weierstrass implications for the localized
polynomial ring $A=K[X]_{(X)}$}

Let $n$ be a positive integer, let $X=\{x_1,\hdots,x_n\}$ be a set
of $n$ variables over a field $K$, and let
$A:=K[x_1,\hdots,x_n]_{(x_1,\hdots,x_n)} = K[X]_{(X)}$ denote the
localized polynomial ring in these $n$ variables over $K$.  Then
the completion of $A$ is $\widehat A=K[[X]]$.

 \proclaim{6.1 Theorem}
 For the localized polynomial ring $A = K[X]_{(X)}$ defined above, if $Q$ is an ideal of
 $\widehat
A$ maximal with  respect to $Q\cap A=(0)$, then $Q$ is a prime
ideal of height $n-1$.
\endproclaim
\demo{Proof} Again it is clear that $Q$ as described in the
statement is a prime ideal. Also the assertion holds  for $n =1$.
Thus  we assume $n \ge 2$. By Proposition~5.4, there exists a
nonzero prime $\p$ in $K[[x_1, x_2]]$ such that $\p \cap K[x_1,
x_2]_{(x_1, x_2)} = (0)$. It follows that $\p\widehat A \cap A =
(0)$. Thus the generic formal fiber of $A$ is nonzero.

Let $P\in \Spec \widehat A$ be a nonzero prime ideal with $P \cap
A = (0)$ and  $\hgt P =r<n-1$.  We construct $Q\in \Spec \widehat
A$ of height $n-1$ with $P\subseteq Q$ and $Q\cap A=(0)$. By
Proposition 2.2, there exists a change of variables $x_1 \mapsto
z_1, \ldots, x_n \mapsto z_n$ and polynomials
$$\alignat2
f_1&\in
K[[z_1,\hdots,z_{n-r}]]\,[z_{n-r+1},\hdots,z_n]\quad&&\text{monic in
$z_n$}\\
f_2&\in K[[z_1,\hdots,z_{n-r}]]\,[z_{n-r+1},\hdots,z_{n-1}]
\quad&&\text{monic
in $z_{n-1}$, etc}\\
 &\vdots && \\
f_r&\in K[[z_1,\hdots,z_{n-r}]]\,[z_{n-r+1}] \quad&&\text{monic in
$z_{n-r+1}$,}
\endalignat
$$
so that $P$ is a minimal prime of $(f_1,\hdots,f_r)\widehat A$ and
$P$ is extended from
$$ R := K[[z_1,\hdots,z_{n-r}]]\,[z_{n-r+1},\hdots,z_n].$$
Let $P_0 := P \cap R$ and extend $f_1,\hdots,f_r$ to a system of
generators of $P_0$, say:
$$P_0 =
(f_1,\hdots,f_r,g_1,\hdots, g_s)R.$$ Using an argument similar to
that  in the proof of Theorem 2.3, write
$$ f_j = \sum_{(i)\in {\Bbb N}^{n-1}} a_{j,(i)} z_2^{i_2}\hdots
z_n^{i_n}\quad \text{ and }\quad  g_k = \sum_{(i)\in {\Bbb
N}^{n-1}} b_{k,(i)} z_2^{i_2}\hdots z_n^{i_n},$$ where $a_{j,(i)},
b_{k,(i)}\in K[[z_1]]$. Let
$$ V_0 := K(z_1,a_{j,(i)}, b_{k,(i)})\cap K[[z_1]].$$
Then $V_0$ is a discrete rank-one valuation domain with completion
$K[[z_1]]$, and $K((z_1))$ has uncountable transcendence degree
over the field of fractions $\Cal Q(V_0)$ of $V_0$.  Let
$\gamma_3,\hdots,\gamma_{n-r}\in K[[z_1]]$ be algebraically
independent over $\Cal Q(V_0)$ and define
$$
\q := (z_3
-\gamma_3z_2,z_4-\gamma_4z_2,\hdots,z_{n-r}-\gamma_{n-r}z_2)K[[z_1,
\ldots, z_{n-r}]].
$$
We see that $\q \cap V_0[[z_2, \ldots, z_{n-r}]] = (0)$ by an
argument similar to that  in \cite{7} and in Claim 2.3.1. Let $R_1
:= V_0[[z_2, \ldots, z_{n-r}]]\,[z_{n-r+1}, \ldots, z_n]$, let $P_1
:= P \cap R_1$ and consider the commutative diagram:

$$
\CD
K[[z_1,\hdots,z_{n-r}]] &\longrightarrow &R/P_0\\
 \uparrow & &\uparrow \\
V_0[[z_2,\hdots,z_{n-r}]] &\longrightarrow &R_1/P_1
\endCD
$$
The horizontal maps are injective finite integral extensions. Let
$W$ be a minimal prime of $(\q, P)\widehat A$. Then $\hgt W = n-2$
and $\q \cap V_0[[z_2, \ldots, z_{n-r}]] = (0)$ implies that $W
\cap R_1 = P_1$. We have found a prime ideal $W \in \Spec \widehat
A$ such that $\hgt W = n-2$, $W \cap A = (0)$ and $P \subseteq W$.
Since $f_1,\hdots,f_r\in W$ and since $\widehat
A=K[[z_1,\hdots,z_n]]$ is the $(f_1,\hdots,f_r)$-adic completion
of $K[[z_1,\hdots,z_{n-r}]]\,[z_{n-r+1},\hdots,z_n]$, the prime
ideal $W$ is extended from
$K[[z_1,\hdots,z_{n-r}]]\,[z_{n-r+1},\hdots,z_n]$.

We claim that $W$ is actually extended from
$K[[z_1,z_2]]\,[z_3,\hdots,z_n]$. To see this let $g\in W\cap
K[[z_1,\hdots,z_{n-r}]]\,[z_{n-r+1},\hdots,z_n]$ and write:
$$ g=\sum_{(i)} a_{(i)} z_{n-r+1}^{i_{n-r+1}}\hdots z_n^{i_n}\in
K[[z_1,\hdots,z_{n-r}]]\,[z_{n-r+1},\hdots,z_n],$$ where the sum is
over all $(i)=(i_{n-r},\dots,i_n)$ and $a_{(i)}\in
K[[z_1,\hdots,z_{n-r}]]$. For all $a_{(i)}$ by Weierstrass we can
write
$$ a_{(i)} = (z_{n-r} -\gamma_{n-r}z_2)h_{(i)} + q_{(i)},$$
where $h_{(i)} \in K[[z_1, \ldots z_{n-r}]]$ and $q_{(i)}\in
K[[z_1,\hdots,z_{n-r-1}]]$.  If $n-r > 3$, we write

$$
q_{(i)} = (z_{n-r-1} - \gamma_{n-r-1}z_2)h'_{(i)} + q'_{(i)},
$$
where $h'_{(i)} \in K[[z_1, \ldots z_{n-r-1}]]$ and $q'_{(i)}\in
K[[z_1,\hdots,z_{n-r-2}]]$. In this way we replace a generating
set for $W$  in $K[[z_1,\hdots,z_{n-r}]]\,[z_{n-r+1},\hdots,z_n]$
by a generating set for  $W$ in $K[[z_1,z_2]]\,[z_3,\hdots,z_n]$.

In particular, we can replace the elements $f_1,\hdots,f_r$ by
elements:
$$\alignat2
h_1&\in K[[z_1,z_2]]\,[z_3,\hdots,z_n]\quad&&\text{monic in
$z_n$}\\
h_2&\in K[[z_1,z_2]]\,[z_3,\hdots,z_{n-1}] \quad&&\text{monic in
$z_{n-1}$, etc}\\
 &\vdots && \\
h_r&\in K[[z_1,z_2]]\,[z_3\hdots,z_{n-r+1}] \quad&&\text{monic in
$z_{n-r+1}$}
\endalignat
$$
and set $h_{r+1} =  z_3-\gamma_3z_2,\hdots, h_{n-2} =
z_{n-r}-\gamma_{n-r}z_2$,  and then extend to a  generating set
$h_1, \ldots, h_{n+s-2}$ for

$$W_0 = W\cap K[[z_1, z_2]]\,[z_3, \ldots, z_n]$$
such that $W_0\widehat A = W$.  Consider the coefficients in
$K[[z_1]]$ of the $h_j$:
$$ h_j = \sum_{(i)} c_{j(i)} z_2^{i_2}\hdots z_n^{i_n}$$
with $c_{j(i)}\in K[[z_1]]$. The set $\{c_{j(i)}\}$ is countable.
Define
$$ V := \Cal Q(V_0)(\{c_{j(i)}\})\cap K[[z_1]]$$
Then $V$ is a rank-one discrete valuation domain  that is
countably generated over $K[z_1]_{(z_1)}$ and $W$ is extended from
$V[[z_2]]\,[z_3, \ldots, z_n]$.

We may also write each $h_i$ as a polynomial in $z_3,\hdots,z_n$
with coefficients in $V[[z_2]]$:
$$ h= \sum \omega_{(i)} z_3^{i_3}\hdots z_n^{i_n}$$
with $\omega_{(i)}\in V[[z_2]]\subseteq K[[z_1, z_2]]$. By the
result of Valabraga \cite{11}, the integral domain
$$ D := \Cal Q(V)(z_2, \{\omega_{(i)}\})\cap K[[z_1, z_2]]$$
is a two-dimensional regular local domain with completion
$\widehat D = K[[z_1, z_2]]$.  Let $W_1 := W \cap D[z_3, \ldots,
z_n]$. Then $W_1 \widehat A = W$. We have shown in Section 5 that
there exists a prime element $q\in K[[z_1,z_2]]$ with $qK[[z_1,
z_2]] \cap D = (0)$. Consider the finite extension
$$ D\longrightarrow D[z_3, \ldots, z_n]/W_1.$$
Let $Q\in \Spec  \widehat A$ be a minimal prime of $(q,W)\widehat
A$. Since $\hgt W = n-2$ and $q \not\in W$,  $\hgt Q = n-1$.
Moreover, $P \subseteq W$ implies  $P\subseteq Q$. We claim that
$$ Q\cap  D[z_3, \ldots, z_n] = W_1 \quad \text{and therefore}\quad Q\cap
A = (0).$$ To see this consider the commutative diagram:

$$
\CD
K[[z_1, z_2]] &\longrightarrow & \quad K[[z_1, \ldots, z_n]]/W \\
 \uparrow & &\uparrow \\
D &\longrightarrow & \quad D[z_3, \ldots, z_n]/W_1\,,
\endCD
$$
which has injective finite horizontal maps. Since $qK[[z_1, z_2]]
\cap D = (0)$, it follows that $Q\cap D[z_3, \ldots, z_n] =  W_1$.
This completes the proof of Theorem 6.1. \qed
\enddemo

\subheading{7. Generic fibers of power series ring extensions}

In this section we apply the Weierstrass machinery from Section 2
to the generic fiber rings of power series extensions.

\proclaim{7.1 Theorem} Let $n \ge 2$ be an integer and let  $y,
x_1, \ldots, x_n$ be variables over the field  $K$. Let $X =
\{x_1, \ldots, x_n \}$. Consider the formal power series ring $R_1
= K[[X]]$  and the extension $R_1 \hookrightarrow
 R_1[[y]] = R$.  Let $U = R_1 \setminus (0)$. For
$P \in \Spec R$ such that $P \cap U = \emptyset$ we have: \roster
\item If $P \not\subseteq XR$, then $\dim R/P = n$ and $P$ is
maximal in the generic fiber  $U^{-1}R$. \item If $P \subseteq
XR$, then there exists $Q \in \Spec R$ such that $P \subseteq Q$,
$\dim R/Q  = 2$ and $Q $ is maximal in the generic fiber
$U^{-1}R$.
\endroster
If $n>2$ for each prime ideal $Q$  maximal in the generic fiber
$U^{-1}R$, we have
$$
\dim R/Q = \cases  n  & \text{ and } R_1 \hookrightarrow R/Q \text{ is finite, or  } \\
2  & \text{ and } Q \subset XR.
\endcases
$$
\endproclaim

\demo{Proof}  Let $P \in \Spec R$ be such that $P \cap U =
\emptyset$ or equivalently $P \cap R_1 = (0)$. Then $R_1$ embeds
in $R/P$. If $\dim(R/P) \le 1$, then the maximal ideal of $R_1$
generates an ideal primary for the maximal ideal of $R/P$. By
\cite{8, Theorem 8.4} $R/P$ is finite over $R_1$, and so $\dim R_1
= \dim(R/P)$, a contradiction. Thus $\dim(R/P) \ge 2.$

If $P \not\subseteq XR$, then there exists a prime element $f \in
P$ that contains a term $y^s$ for some positive integer $s$. By
Weierstrass, it follows that $f = g\epsilon$, where $g \in
K[[X]]\,[y]$ is a nonzero monic polynomial in $y$ and $\epsilon$ is
a unit of $R$. We have $fR = gR \subseteq P$ is a prime ideal and
$R_1 \hookrightarrow R/gR$ is a finite integral extension. Since
$P \cap R_1 = (0)$, we must have $gR = P$.

If $P \subseteq XR$ and $\dim(R/P)> 2$, then Theorem 2.3
 implies there exists $Q \in \Spec
R$ such that $\dim(R/Q) = 2$, $P \subset Q \subset XR$ and $P \cap
R_1 = (0) = Q \cap R_1$, and so $P$ is not maximal in the generic
fiber. Thus $Q \in \Spec R$ maximal in the generic fiber of $R_1
\hookrightarrow R$ implies that the dimension of $\dim(R/Q)$ is
$2$, or equivalently that $\hgt Q=n-1$. \qed
\enddemo

\proclaim{7.2 Theorem} Let $n$ and $m$ be positive  integers,  and
let $X = \{x_1, \ldots, x_n\}$ and $Y = \{y_1, \ldots, y_m\}$ be
sets of independent variables over the field  $K$. Consider the
formal power series rings $R_1 = K[[X]]$ and  $R = K[[X, Y]]$ and
the extension $R_1 \hookrightarrow
 R_1[[Y]] = R$. Let $U = R_1 \setminus (0)$. Let $Q \in \Spec R$  be maximal with
respect to $Q \cap U = \emptyset$. If $n = 1$, then $\dim R/Q = 1$
and $R_1 \hookrightarrow R/Q$ is finite.

If $n \ge 2$, there are two possibilities \roster \item $R_1
\hookrightarrow R/Q$ is finite, in which case $\dim R/Q = \dim R_1
= n$, or \item $\dim R/Q = 2$.
\endroster

\endproclaim

\demo{Proof} First assume $n = 1$, and let $x = x_1$. Since $Q$ is
maximal with respect to $Q \cap U = \emptyset$, for each $P \in
\Spec R$ with $Q \subsetneq P$ we have $P \cap U$ is nonempty and
therefore $x \in P$. It follows that $\dim R/Q = 1$, for
otherwise,
$$
Q = \bigcap\{ P\,\, |\,\, P \in \Spec R \text{ and } Q \subsetneq
P \, \},
$$
which implies $x \in Q$. By \cite{8, Theorem 8.4}, $R_1
\hookrightarrow R/Q$ is finite.

It remains to consider the case where  $n \ge 2$. We proceed by
induction on $m$. Theorem 7.1 yields the assertion for $m = 1$.
Suppose $Q \in \Spec R$ is maximal with respect to $Q \cap U =
\emptyset$.  As in the proof of Theorem 7.1, we have $\dim R/Q \ge
2$. If $Q \subseteq (X, y_1, \ldots, y_{m-1})R$, then by Theorem
2.3 with $R_0 = K[y_m]_{(y_m)}[[X, y_1, \ldots, y_{m-1}]]$, there
exists $Q' \in \Spec R$ with $Q \subseteq Q'$, $\dim R/Q' = 2$,
and $Q \cap R_0 = Q' \cap R_0$. Since $R_1 \subseteq R_0$, we have
$Q' \cap U = \emptyset$. Since $Q$ is maximal with respect to $Q
\cap U = \emptyset$, we have $Q = Q'$, so $\dim R/Q = 2$.

Otherwise, if   $Q \not\subseteq (X, y_1, \ldots, y_{m-1})R$, then
there exists a prime element $f \in Q$ that contains a term
$y_m^s$ for some positive integer $s$. Let $R_2 = K[[X, y_1,
\ldots,y_{m-1}]]$. By Weierstrass, it follows that $f =
g\epsilon$, where $g \in R_2[y_m]$ is a nonzero monic polynomial
in $y_m$ and $\epsilon$ is a unit of $R$. We have $fR = gR
\subseteq Q$ is a prime ideal and $R_2 \hookrightarrow R/gR$ is a
finite integral extension. Thus $R_2/(Q \cap R_2) \hookrightarrow
R/Q$ is an integral extension. It follows that $Q \cap R_2$ is
maximal in $R_2$ with respect to being disjoint from $U$. By
induction $\dim R_2/(Q \cap R_2)$ is either $n$ or $2$. Since
$R/Q$ is integral over $R_2/(Q \cap R_2)$, $\dim R/Q$ is either
$n$ or $2$. \qed
\enddemo

\subheading {7.3 Remark} In the notation of Theorem 1.1, Theorem
7.2 proves the second part of the theorem, since $\dim R=n+m$.
Thus if $n=1$, $\hgt Q=m$. If $n\ge 2$, the two cases are (i)
$\hgt Q=m$ and (ii) $\hgt Q=n+m-2$, as in (a) and (b) of Theorem
1.1.

\medskip

Using the TGF terminology discussed in the introduction, we have
the following corollary to Theorem 7.2.

\proclaim{7.4 Corollary} With the notation of Theorem 7.2, assume
$P\in\Spec R$ is such that $R_1 \hookrightarrow R/P =: S$ is a TGF
extension. Then $\dim S = \dim R_1 = n$ or $\dim S = 2$.
\endproclaim

\enddocument